\documentclass[a4paper]{amsart}
\usepackage{amssymb}
\usepackage{verbatim}
\usepackage{amscd}
\def\ent{{{\operatorname{h}}}}
\def\card{{{\operatorname{card}}}}

\def\id{{{\operatorname{id}}}}

\numberwithin{equation}{section}
\theoremstyle{plain}
\newtheorem{theorem}[equation]{Theorem}

\newtheorem{lemma}[equation]{Lemma}
\newtheorem*{(DQ1)}{(DQ1)}
\theoremstyle{definition}

\theoremstyle{remark}

\begin{document}
\title [On  images of sofic systems ]{On  images of sofic systems}
\author{Wolfgang Krieger}
\dedicatory{Dedicated to the memory of Roy Adler}
\begin{abstract}
Let $\Sigma$ and $\bar\Sigma$  be finite alphabets. 
For topologically transitive sofic systems $ X\subset  \Sigma^{\Bbb Z}$ and   
$\widetilde X\subset  \widetilde\Sigma^{\Bbb Z}$  
we give a necessary and sufficient condition for the existence of a homomorphism from $X$ to 
$\widetilde X$. For topologically mixing sofic systems $X \subset  \Sigma^{\Bbb Z}$ and  $\widetilde X\subset  \widetilde\Sigma^{\Bbb Z}$, such that the topological entropy of 
$\widetilde X$ is less than the topological entropy of $X$, we give a necessary and sufficient condition for the existence of a homomorphism of  $X$ onto $\widetilde X$.
\end{abstract}
\maketitle

\section{Introduction}
Let $\Sigma $ be a finite alphabet.
On the shift space $ \Sigma^{\Bbb Z}$ there acts the full shift $S$,
$$
S((x_{i})_{i\in \Bbb Z}) = (x_{i+1})_{i\in \Bbb Z},\quad   (x_{i})_{i\in \Bbb Z} \in \Sigma^{\Bbb Z}.
$$
A closed shift invariant subset $X$ of $\Sigma^{\Bbb Z}$ together with the restriction $S_X$ of the shift $S$ to $X$ is a dynamical system that is called a subshift.  For introductions to the theory of subshifts see \cite {LM} and \cite {Ki}.
By a homomorphism (or factor map) $\varphi: X \to  \widetilde X $ from a subshift 
$X \subset  \Sigma^{\Bbb Z}$ to a subshift $\widetilde X\subset \widetilde\Sigma^{\Bbb Z}$  is meant  a continuous shift-commuting map from $X$  into  $\widetilde X$. 
One of the basic classes of  subshifts are the subshifts of finite type, that are defined as the subshifts that are obtained from a finite set $\mathcal F$ of words by excluding from a shift space $ \Sigma^{\Bbb Z}$ the words in $\mathcal F$ \cite {P}. The class of subshifts of finte type is closed under topological conjugacy. Sofic systems \cite {W} are the homomorphic images of subshifts of finite type. 

The study of topological conjugacy and of related notions of equivalence, together with the study of homomorphisms, 
has continuously been a central topic in this area of symbolic dynamics since its beginnings \cite{P, PT}(among  the references in \cite {PT} see \cite{AGW, AM}, and see also \cite {AKM}).
In this paper we are concerned with homomorphisms  between 
topologically transitive sofic systems and surjective homomorphisms between 
topologically mixing sofic systems.
For a sofic system $X \subset  \Sigma^{\Bbb Z}$ and  a topologically transitive aperiodic subshift of finite type 
$ \widetilde X \subset  \widetilde \Sigma^{\Bbb Z}$, such that the topological entropy of  $X$ exceeds the topological entropy of $\widetilde X$, Mike Boyle  has shown that the periodic point condition is sufficient for the existence of a surjective 
homomorphism of  $X$ onto $ \widetilde X$ \cite [Corollary (2.6)]{B}.  The periodic point condition says, that every $\pi\in \Bbb N$ that appears as the period of a periodic point of $X$ has a divisor that appears as a period of a periodic point of $\widetilde X$. Boyle also gave a sufficient condition for the existence of a surjective homomorphism of a topologically transitive sofic system onto a topologically mixing sofic system of lower entropy 
 \cite [Theorem  3.3]{B}. This  was strengthened  by Klaus Thomsen  \cite [Theorem 9.13]{T}  and further strengthened  by Jan Nielsen  \cite {N}.

In section 4 of this paper we give for topologically transitive sofic systems 
$X \subset  \Sigma^{\Bbb Z}$ and  $\widetilde X\subset \widetilde\Sigma^{\Bbb Z}$ 
a necessary and sufficient condition for the existence of a homomorphism from $X$ to 
$\widetilde X$. 
In Section 5  we point out, that this condition is decidable. This is not unexpected, since the basic properties of regular languages are decidable \cite [Chapter 3]{HU}, and also in view of the results on algorithms for sofic systems \cite {CP}.

In Section 6 we will obtain for  topologically mixing sofic systems $X \subset  \Sigma^{\Bbb Z}$ and  $\widetilde X\subset \widetilde\Sigma^{\Bbb Z}$ , such that the topological entropy of $X$ exceeds the topological entropy of $\widetilde X$, a necessary and sufficient condition for the existence of a surjective homomorphism of  $X$ onto $\widetilde X$. We point out that this condition is also decidable.

As in  \cite{Kr} our arguments are based on the construction of a compact-open set such that the entries of any point into the set under the action of the shift are sufficiently spaced. This set is used to separate periodic events from non-periodic events. Also, we create  sufficiently long periodic events  by the use of certain auxiliary maps, that we introduce in Section 3. In a preliminary Section 2 we introduce notation and recall some tools. At certain stages of our constructions we find it necessary to choose a time direction. This choice is always arbitrary.

Acknowledgement: Thanks go to Jan Nielsen for a discussion on the results of his thesis,  and to Mike Boyle for his assistance. 

\section{Preliminaries}
 We  introduce  notation and terminology for subshifts $ X\subset  \Sigma^{\Bbb Z}$. We set
$$
x_{[i,k]} = (x_{j})_{ i \leq j \leq k}, \qquad  i, k \in \Bbb Z, j \leq k,\qquad
 (x_{i})_{i \in\Bbb Z} \in X,
$$
and  we set
$$
X_{[i,k]}= \{ x_{[i,k]} : x \in X \},  \qquad  i, k \in \Bbb Z, j \leq k,\qquad
 (x_{i})_{i \in\Bbb Z} \in X,
$$
and we use  similar notation in the case that indices range over open, half-open or semi-infinite intervals. The set of admissible words of a subshift $X \subset \Sigma^{\Bbb Z}$ we denote
by $\mathcal L (X)$. 
The topological entropy of a subshift is denoted by $\ent$.
The symbol  $\ell$ denotes the length of a word. The empty word we denote by $\epsilon$. For cylinder sets we use the notation like
$$
Z(b) =\{ x \in X: x_{[1, \ell(b)]} = b \}, \quad b \in \mathcal L(X).
$$

We set 
$$
\Gamma^-(w)= \{ u \in \mathcal L(X): uw \in  \mathcal L(X) \}, \qquad w \in \mathcal L(X).
$$
The meaning of $ \Gamma^+$ is symmetric. We set
$$
\Gamma(w)= \{(u, v) \in  \mathcal L(X)^2: uwv \in \mathcal L(X) \}, 
\qquad w \in \mathcal L(X).
$$
The set $\Gamma(w)$ is called the context of $w$.  We say that $ w\in \mathcal L(X)$ and 
$ w^{\prime}\in \mathcal L(X)$  belong to the same context class if 
$\Gamma(w)= \Gamma(w^\prime)$ and we write $w \approx w^{\prime}$.
The set $[\mathcal L(X)]_\approx$ with a zero adjoined (excepting the case of a full shift)
with the product given by
$$
[u]_\approx[v] _\approx=
\begin{cases}
[uv]_\approx, &\text{if  $uv \in \mathcal L(X)$},  \\
0, &\text {if $uv \not\in \mathcal L(X)$},
\end{cases}
$$
is a semigroup (the syntactic semigroup of $X$). Sofic systems have a finite syntactic semigroup.
We denote by $\Delta(X)$ the set of $\delta \in [ \mathcal L(X)]_\approx$ such that there exists a word
  $c\in \mathcal L(X)$, such that
\begin{align*} 
 \delta[c]_\approx =   \delta. \tag {2.1}
\end{align*}
We order the set $\Delta(X)$ linearly.

A word $u \in\mathcal L(X)$ is called synchronizing if $\Gamma (u) =  \Gamma_- (u) \times\Gamma_+ (u) $ 
We denote the set of context classes of synchronizing words of $X$ by $\Delta_{synchro}(X)$. We note that for topologically transitive subshifts $X$
$$
\Delta_{synchro}(X)\subset \Delta(X).
$$

We denote the set of periodic points of $X$ by $P(X)$, and we denote the smallest period of $p \in  P(X)$ by $\pi(p)$. We set
$$
P_{\langle K \rangle}(X)= \{ p \in P(X) : \pi(p) \leq k \}, \quad k \in \Bbb N.
$$

We set
$$
\mathcal A(X) = \{ p_{[0, \pi(p))} : p \in P(X)\},
$$
$$
\mathcal A_{\langle K \rangle}(X)= \{ a \in \mathcal P(X) : \ell(a) \leq k \}, \quad k \in \Bbb N.
$$
We say that $a, a^{\prime}\in \mathcal A(X)$ are conjugate and write 
$a \sim  \thinspace a^{\prime}$ if $a^{\prime}$ is obtained from 
$a$ by cyclically permuting the symbols of $a$. 
Given a word $a \in \mathcal A(X)$ we use the notation 
$p^{(-, a)}$ ( $p^{(+, a)}$ ) for the
element in $X_{(- \infty. 0]}$ ( $X_{[1, \infty )}$)
that carries the semi-infinite concatenation of $a$.
For $a \in \mathcal A(X)$ and $k\in \Bbb N$ we set
$$
\tau ^{(k)}(a) = p^{(+, a)} _{(k, k + \ell(a))}.
$$
For $a\in \mathcal A(X)$ and $u \in \mathcal L(Y)$ we set
$$
\kappa_a(u) = \max \{I \in [1, \ell(u)]: u_{[1, I]} = p^{(+, a)}_{[1, I]} \}.
$$

For $a \in \mathcal A(X)$ we denote by $\gamma (a)$ the smallest $N \in \Bbb N$ such that there is a $K \in \Bbb N$, such that
$$
[a^{K + N}]_\approx = [a^{K}]_\approx. 
$$
\begin{lemma}
For conjugate  $a, \bar a \in \mathcal A(X)$ one has that $\gamma (a) = \gamma (\bar a).$ 
\end{lemma}
\begin{proof}
Let  $N \in \Bbb N$ be such that
$
[a^{ N}]_\approx = [a^{N + \gamma(a)}]_\approx, 
$
and let $r$ be given by
$
\bar a= \tau^{(r)}(a).
$
Then
$$
[ \bar a ^{N + 1 + \gamma(a)}]_\approx = [\bar a_{[1, r]}a ^{N  + \gamma(a)}
\bar a_{(r, \ell(a)]}]_\approx=
 [\bar a_{[1, r]}a ^{N}
\bar a_{(r, \ell(a)]}]_\approx = [ \bar a ^{N + 1 }]_\approx. \qed
$$
\renewcommand{\qedsymbol}{}
\end{proof}
\noindent
We set
$$
\mu(a) = \gamma (a)\  \max_{a_\circ \sim a}\left( \min \thinspace \{N \in \Bbb N:  [a_\circ^{ N\gamma(a)}]_\approx 
= [a_\circ^{( N+1)\gamma(a)}]_\approx \}\right),
\qquad a \in \mathcal A(X).
$$

Given a subshift $X \subset \Sigma^{\Bbb Z}$ and $H \in  \Bbb N$, we denote for $h\in [1,H]$ by 
$\mathcal P_{2H, h}(X)$ the set of words 
$w \in \mathcal L_{2H}(X)$, such that
$$
w_i = w_{i-h}, \qquad h \leq i \leq 2H,
$$
and set
$$
\mathcal P_{2H}(X) = \bigcup_{1 < h \leq H}\mathcal P_{2H, h}(X).
$$
For $w \in\mathcal P_{2H}(X)$ we denote by  $\pi(w)$ the smallest $h \in [1, 2H]$, such that 
$
w \in \mathcal P_{2H, h}(X)
$
and we set 
$$
g(w) = w_{[1,\pi(w) ]},  \qquad w \in \mathcal P_{2H}(X).
$$

By  \cite [Lemma 2]{Kr} that, given an $H \in \Bbb N$ there is a set $F\subset X$, such that 
\begin{align*}
F \cap S_X^h (F)= \emptyset, \qquad 1 < h <H, \tag {2.2}
\end{align*}
and such that
\begin{align*}
\bigcup_{-H < h < H}S_X^h (F) \supset 
 \bigcup_{w \in \mathcal L_{2H}(X)\setminus \mathcal P_{2H}(X)}  \tag {2.3}
Z(w).
\end{align*}

For a compact-open set $A \subset X$ and for $x\in X$ we introduce the notation (in which we suppress the symbol $X$):
$$
i_A^+(x) =
\begin{cases}
i, &\text{if  $x \in (\bigcap_{1\leq j <i} (X - S_X^j(A)))\cap S_X^i(A)$},  \\
\infty, &\text {if $x \in \bigcap_{1\leq j <\infty} (X - S_X^j(A)) $}, \qquad \qquad\qquad x \in A.
\end{cases}
$$
The meaning of $i_A^-(x)$ is symmetric. 

Given  subshifts $ X\subset  \Sigma^{\Bbb Z}$ and   
$\widetilde X\subset  \widetilde\Sigma^{\Bbb Z}$ and a homomorphism $\varphi: X \to \widetilde X$, we set 
$$
^\varphi a = \varphi (p^{(-, a)} , p^{(+, a)})_{[1, \pi(\varphi(p^{-, a} , p^{(+, a)}))]}
 \qquad a\in \mathcal A(X).
$$
Also we recall that, given a homomorphism 
 $\varphi:X \to \bar{X}$, 
 there is for some $L \in \Bbb Z_+$ a block map
 $$\bold\Phi:X_{[-L,L]} \to \widetilde{\Sigma},$$
 that implements the homomorphism $\varphi$ with the action of $\bold\Phi$ on words given by
 $$
\varphi(x)_{[I_-+L,I_+-L]} = (\bold\Phi(x_{[i-L,i+L]}))_{{ I_-+L\leq i \leq I_+-L}}, \ \   x \in X,  I_-, I_+ \in \Bbb Z, I_+ -I_- > 2L.
$$

\section{Auxiliary maps}
We consider a topologically transitive sofic system $ X \subset \Sigma^{\Bbb Z}$. 
We denote for  $\delta \in \Delta(X)$ by $\eta(\delta)$ the shortest length of a $c \in \mathcal L(X)$ such that (2.1) holds.
We also select for  $\delta \in \Delta(X)$ a word $c(\delta)$ of length $\eta(\delta)$ such that (2.1) holds. We have then automatically also selected for $\delta \in \Delta(X)$ an 
$f(\delta) \in \mathcal A (X)$ which is given by
$$
c(\delta) = f(\delta)^{\frac{\eta(\delta)}{\ell(f(\delta))}}.
$$
Also denote for $\delta \in \Delta(X)$ by $\mathcal M(\delta)$  the set of $m \in \Bbb N$, such that there is a 
$c \in \mathcal L _m(X)$ such that (2.1) holds, and denote by $\nu(\delta)$ the smallest $n\in \Bbb N$ such that
$$
\mathcal M (\delta)  \cap [n\eta(\delta), (n + 1)\eta(\delta) ] = \mathcal M (\delta)\mod \eta(\delta).
$$

Set
$$
\mathcal U(X) =\{u \in \mathcal L(X): \ell (u) > \card (\Delta(X))\},
$$
and denote for  $u \in \mathcal U(X) $ the smallest index $I \in [1, \card (\Delta(X))]$, such that
$$
[U_{[1,I]}]_\approx \in \Delta(X) ,
$$
by $I^-(u)$, and set
$$
\delta (u) = [u_{[1,I^-(u)]}]_\approx.
$$
Set
$$
I^+(u) = \max \{I \in ( I^-, \ell (u)]: 
u_{(I^-(u), I]} = p_{[ 1, I -   I^-(u)]}^{(+, f(\delta(u) ) )} \}, \qquad u \in \mathcal U(X).
$$
For $\delta \in\Delta(X)$ we denote by  $\mathcal U_\delta(X) $ the set of $u \in \mathcal U$ such that there exist indices
$$
J \in  (I^+(u) , \ell (u)  ],
$$
such that
$$
\delta = [u_{[I^+(u), J]}]_\approx,
$$
and such that one has for the smallest such index $J^-(u, \delta)$ and the largest such index
 $J^+(u, \delta)$, that
$$
J^+(u, \delta) - J^-(u, \delta) \geq \eta(\delta)(\nu(\delta)  + \eta(\delta(u)).
$$
We define maps
$$
\Psi_\delta: \mathcal U_\delta(X) \to \mathcal U(X), \qquad \delta \in \Delta,
$$
that leave the length and the context class of a word unchanged. For this purpose we select for every 
$\delta \in \Delta$ and every $m \in \mathcal M(\delta) $,
such that
$$
m\geq \nu(\delta)\eta(\delta),
$$
 a word $b_m(\delta) \in \mathcal M_m(X)$ such that
 $$
 \delta[b_m(\delta)]_\sim = \delta.
 $$
 From the word $u\in  \mathcal U_\delta(X)$ we obtain the word $\Psi_\delta(u)$ by inserting in $u$ after  index $I^-(u)$ the word $c(\delta(u))^{\eta(\delta)}$, and by replacing in $u$ the word
$
u_{ (J^-(u, \delta), J^+(u, \delta) ]}
$ by the word
$$
b_{J^+(u, \delta) - J^-(u, \delta)- \eta(\delta)\eta(\delta(u))}(\delta).
$$

We define a map
$$
\Psi: \mathcal U (X) \to \mathcal U(X) \setminus \bigcup_{\delta\in \Delta}\mathcal U_\delta(X),
$$
that leaves length and context class of a word unchanged.
We set
$$
\Psi \restriction  \mathcal U(X) \setminus \bigcup_{\delta\in \Delta}\mathcal U_\delta(X) = \id
$$
and otherwise obtain for
$
u \in \bigcup_{\delta\in \Delta}\mathcal U_\delta(X),
$
its image under $\Psi $
with its specific properties  by an inductive procedure. 
For 
$u\in \bigcup_{\delta\in \Delta}\mathcal U_\delta(X)$
we denote by $\delta_\circ(u)$ the smallest $\delta\in \Delta$ such that $u\in \mathcal U_\delta$.
We determine for
$$
u \in  \bigcup_{\delta\in \Delta}\mathcal U_\delta(X).
$$
a $K(u)\in [1, \ell(u) - I^-(u)]$,
and 
$$
 u_k\in \mathcal U(X), \delta^{(k)}\in \Delta(X),  \qquad 1 \leq k \leq K(u),
 $$
 by
\begin{align*}
&u_1 = u, \delta^{(1)}= \delta_\circ(u),
\\
&u_k = \Psi_{ \delta^{(k-1)}}(u_{k-1}), 
 \delta^{(k)} = \delta_\circ(u_{k-1} ),  \quad 1 < k <  K(u) ,
\\
&u_{K(u)-1} \in \bigcup_{\delta \in \Delta(X)} \mathcal U_{ \delta}(X),
\\
&u_{K(u)} = \Psi_{\delta^{(K(u)-1)}}(u_{K(u)-1}) \in 
 \mathcal U(X) \setminus \bigcup_{\delta\in \Delta}\mathcal U_\delta(X),  
\end{align*}
and we set
$$
\Psi(u) = u_{K(u)}.
$$
\begin{lemma}
For $\delta \in \Delta(X)$ and for 
\begin{align*}
u\in  \mathcal U(X) \setminus \bigcup_{\delta\in \Delta}\mathcal U_\delta(X)
\tag{3.1}
\end{align*}
one has
$$
\ell(u) -I^+(u)   \leq \card ( \Delta(X)) (1 +
\max_{\delta \in \Delta}\{\eta(\delta)(\nu(\delta)  + \eta(\delta(u))\}).
$$
\end{lemma}
\begin{proof} 
If there are indices $J_-, J_+ $,
$$
I^+(u) < J_- < J_+ < \ell(u),
$$
such that
$$
J_+ - J_- \geq \max_{\delta \in \Delta}\{\eta(\delta)(\nu(\delta)  + \eta(\delta(u))\} - 1,
$$
and
$$
\{[u_{[I^+(u),j]}] _\approx: 1 \leq j \leq J_-\} \cap \{[u_{[I^+(u),j]}] _\approx : J_+ \leq j \leq \ell (a)\} \neq \emptyset,
$$
then (3.1) cannot hold.
\end{proof}

We also define maps
$$
\Psi_a: \Gamma^+(a^{\mu(a)}) \to \Gamma^+(a^{\mu(a)}) , \qquad a \in \mathcal A(X),
$$
that leave length and context class of a word unchanged. The construction of the maps  
$\Psi_a$ is very similar to the construction of the map $\Psi$. We denote for $\delta \in \Delta$ and $a \in \mathcal A(X)$ by 
$\mathcal U_{\delta, a}$ the set of $u\in \Gamma^+(a^{\mu(a)})$, such that there are indices 
$J\in (\kappa_a(u), \ell (u))$, such that  $\delta = [u_{[1,J]}]_\approx$, and such that one has for the smallest such index  $J^-(u, \delta)$ and the largest such index $J^+(u, \delta)$, that 
$$
J^+(u, \delta) - J^-(u, \delta) \geq \eta(\delta)(\nu(\delta )+ \ell(a)\gamma(a)).
$$
We denote by $\Psi_{a,\delta}(u)$ the word that is obtained from the word 
$u\in \Gamma^+(a^{\mu(a)})$ by appending the word 
$a^{\gamma(a) \eta(\delta(u))}$ on the left and by replacing in  the word $u$ the word  
$$
u_{( J^-(u, \delta) , J^+(u, \delta) ]}
$$
 by the word
$$
b (J^+(u, \delta) - J^-(u, \delta)- \ell(a)\gamma(a) \eta(\delta(u))).
$$
For 
$u\in \bigcup_{\delta\in \Delta}\mathcal U_{a,\delta}(X)$
we denote by $\delta^\circ(u)$ the smallest $\delta\in \Delta$ such that
$u\in \mathcal U_{a,\delta}$.
We set
$$
\Psi_a \restriction  \Gamma^+(a^{\mu(a)})\setminus
 \bigcup_{\delta \in \Delta} \mathcal U_{a, \delta}(X)
= \id,
$$
and for $u \in  \bigcup_{\delta \in \Delta} \mathcal U_{a, \delta}(X) $ we have the inductive procedure
\begin{align*}
&u^{(1)} = u, \delta_{1}= \delta^\circ(u),
\\
&u^{(k)} = \Psi_{ \delta_{k-1}}(u^{(k-1)}), 
 \delta_{k} = \delta^\circ(u_{k-1} ),  \quad 1 < k <  K(u) ,
\\
&u^{(K(u)-1)} \in \bigcup_{\delta \in \Delta(X)} \mathcal U_{a, \delta}(X),
\\
&u^{K(u)} = \Psi_{\delta_{K(u)-1}}(u^{(K(u)-1)}) \in 
 \mathcal U(X) \setminus \bigcup_{\delta\in \Delta}\mathcal U_{a,\delta}(X),  
\end{align*}
and we set
$\Psi_a(u) = u^{K(u)} $.

The following lemma is  analogous to  Lemma (3.1).

\begin{lemma}
For $a\in \mathcal A(X)$ and for 
\begin{align*}
v\in  \Gamma^+(a^{\mu(a)})\setminus\bigcup_{\delta \in \Delta(X)} 
\mathcal U_{a, \delta}(X),
\end{align*}
one has
$$
\ell(v)- \kappa_a(v) \leq \card ( \Delta(X)) 
\max_{\delta \in \Delta}\{\eta(\delta)(\nu (\delta) + \ell(a)\gamma(\delta))\}.
$$
\end{lemma}

\section{Homomorphisms}
We consider a topologically transitive sofic system $X \subset \Sigma^{\Bbb Z}$. We denote by 
$\Omega_{\langle H \rangle}(X)$ the set of triples
$$
(a[-],c,a[+]) \in \mathcal A_{\langle H \rangle}(X) \times (\mathcal L(X) \cup \{\epsilon\}) \times 
\mathcal A_{\langle H \rangle}(X),
$$
such that
$$
a[-]^{\mu(a[-])}ca[+]^{\mu(a[+])}\in \mathcal L(X),
$$
and such that one has, in the case, that $c \neq \epsilon$, that 
$$
c \in \Gamma^+(a[-]^{\mu((a[-])})\setminus\bigcup_{\delta \in \Delta(X)}
 \mathcal U_{a[-], \delta}(X),
$$
and that the initial symbols of $c$ and $a[-]$ are different, and the final symbols of $c$ and 
$a[+]$ are different.  For a given $(a[-],c,a[+])\in \Omega_{\langle H \rangle}(X)$
 we define a point $z^{(a[-],c,a[+])}\in X$ 
 by
$$
z^{(a[-],c,a[+])}_{(-\infty, 0]}=p^{(-,a[-])}c , \quad z^{(a[-],c,a[+])}_{[1, \infty)}
=p^{(+,a[+])}.
$$

\begin{lemma}
Let $x \in X\setminus P(X)$ be left asymptotic to $ p^- \in P(X)$ and right asymptotic to $p^+\in P(X) .$
Then there exist unique $t \in \Bbb Z$ and   $(a[-],c,a[+]) \in \Omega_{\langle H \rangle}(X)$
such that
$$
x = S_X^{-t} (z^{(a[-],c,a[+])}).
$$
\end{lemma}
\begin{proof}
Let $k^+(x) \in \Bbb Z$ be given by the condition
$$
x_{[ k^+(x), \infty)} =p_{[ k^+(x), \infty)} ^+,
\qquad
x_{[ k^+(x)-1, \infty)} \neq p_{[ k^+(x)-1, \infty)} ^+.
$$
In the case that 
\begin{align*}
x_{(-\infty, k^+(x))} =p^- _{(-\infty, k^+(x))} , \tag {4.1}
\end{align*}
$c$ is empty, and
\begin{align*}
&t = k^+(x),  \\
& a[-]= x_{[k^+(x)- \pi(p^-), k^+(x))}, \  \  a[+] =x_{[k^+(x), k^+(x)+\pi(p^+) )}.
\end{align*}
In the case that (4.1) does not hold, let
$k^-(x) \in \Bbb Z$ be given by the condition
$$
\qquad
x_{(-\infty,k^-(x))} = p_{(-\infty,k^-(x))}^-,
\qquad
x_{(-\infty,k^-(x)+1)} \neq p_{  (-\infty,k^-(x)+1)}^-,
$$
and have
\begin{align*}
&t= k^+(x) , \\
&a[-]=  x_{[k^-(x)- \pi(p^-), k^-(x))}  , \
c = x_{[k^-(x), k^+(x)]} , \
a[+] = x_{(k^+(x), k^+(x)+\pi(p^+) )}.  
\qed
\end{align*}
\renewcommand{\qedsymbol}{}
\end{proof}

For a pair
$$
((a[-],c,a[+]) ,(a^\prime[-],c^\prime,a^\prime[+]) ) \in 
(\Omega_{\langle H \rangle}(X) \setminus \{(a, \epsilon, a): a  \in \mathcal A_{\langle H \rangle}(X)\})^2,
$$
such that $ a[+]$ and $ a^\prime[-]$ are conjugate, we denote by 
$$
\rho ((a[-],c,a[+]) ,(a^\prime[-],c^\prime,a^\prime[+]) ) 
$$
the set of $r\in [1, \ell(a[+])]$, such that the initial symbol of $\tau^{(r)}(a[+] )$ is different from the initial symbol of $c^\prime$ (in other words, such that 
$(\tau^{(r)}(a[+] ) ), c^\prime,a^\prime[+]))\in  
\Omega_{\langle H \rangle}(X)$).
We denote by 
$\Pi_{\langle H \rangle}(X)$ the set of pairs 
$$
((a[-],c,a[+]) ,(a^\prime[-],c^\prime,a^\prime[+]) ) \in 
(\Omega_{\langle H \rangle}(X) \setminus \{(a, \epsilon, a): a  \in \mathcal A_{\langle H \rangle}(X)\})^2,
$$
such that
$$
a^\prime[-] = \tau^{(\min(\rho ((a[-],c,a[+]) ,(a^\prime[-],c^\prime,a^\prime[+])))}(a[+]).
$$
Given a triple 
$(a[-],c,a[+])\in \Omega_{\langle H \rangle}(X)$, we set
$$
\beta((a[-],c,a[+])) = gcd(a[-],a[+]).
$$

Given sofic systems  $ X\subset  \Sigma^{\Bbb Z}$ and   
$\widetilde X\subset  \widetilde\Sigma^{\Bbb Z}$, and $H\in \Bbb N$, such that 
$P_{\langle H \rangle}(X)\neq \emptyset$,
we consider maps
$$
\Theta_{\langle H \rangle}(X,\widetilde X): 
\Omega_{\langle H \rangle}(X) \setminus \{(a, \epsilon, a): a  \in \mathcal A_{\langle H \rangle}(X)\}
\to\Omega_{\langle H \rangle}(\widetilde X)\times \Bbb Z_+,
$$
that assign to a triple 
$$
(a[-],c,a[+]) \in
\Omega_{\langle H \rangle}(X) \setminus \{(a, \epsilon, a): a  \in \mathcal A_{\langle H \rangle}(X)\}
$$
an element 
$
((\widetilde a[-],\widetilde c,\widetilde a[+]),\widetilde t)
$
of
$
\Omega_{\langle H \rangle}(\widetilde X)\times \Bbb Z_+,
$
such that 
$$
0 \leq \widetilde t < \beta(((\widetilde a[-],\widetilde c,\widetilde a[+])).
$$
Given a map $\Theta_{\langle H \rangle}(X,\widetilde X)$  we denote for
$$
((a[-],c,a[+]) ,(a^\prime[-],c^\prime,a^\prime[+]) ) \in \Pi_{\langle H \rangle}(X)
$$
the image under $\Theta_{\langle H \rangle}(X,\widetilde X)$  of $(a[-],c,a[+])$ by 
$((\widetilde a[-],\widetilde c,\widetilde a[+]), \widetilde t)$, and the image under 
$\Theta_{\langle H \rangle}(X,\widetilde X)$  of $(a^\prime[-],c^\prime,a^\prime[+])$ by 
$(\widetilde a^\prime[-],\widetilde c^\prime,\widetilde a^\prime[+]), \widetilde t^\prime)$,
and for
$$
r \in \rho ((a[-],c,a[+]) ,(a^\prime[-],c^\prime,a^\prime[+]) )
$$
we denote by 
$$
 \Xi
({\Theta_{\langle H \rangle}(X,\widetilde X)},(a[-],c,a[+]),(a^\prime[-],c^\prime,a^\prime[+]), r)$$
 the set of remainders that are left by the numbers in
$$
\{r + \ell(c^\prime) - \ell(\widetilde c^\prime) + \widetilde t^\prime - \widetilde t + m: 
0\leq m <\gamma(\widetilde a[+])\ell(\widetilde a[+])\}
$$
after division by
$
\gamma(\widetilde a[+])\ell(\widetilde a[+]).
$

Given $H \in \Bbb N$, such that $ P_{\langle H \rangle}(X)\neq \emptyset$, and a homomorphism
$$
\varphi_\circ: P_{\langle H \rangle} (X) \to P_{\langle H \rangle}( \widetilde X),
$$
we say, that the map $\Theta_{\langle H \rangle}(X,\widetilde X)$ accompanies $\varphi_\circ$, if the following conditions 
$(A)$ and  $(B)$ are satisfied:

\medskip

\noindent
$(A)$ For 
$(a[-],c,a[+]) \in \Omega_{\langle H \rangle}(X)$ one has
$\widetilde a[-] = \tau^{(\widetilde t)} (a[-]),$ and $ \widetilde a[+] = 
\tau^{(\widetilde t)} (a[+])$.

\medskip

\noindent
$(B)$ For
$$
(a^{(k)}[-],c^{(k)},a^{(k)}[+]) \in
\Omega_{\langle H \rangle}(X) \setminus \{(a, \epsilon, a): a  \in \mathcal A_{\langle H \rangle}(X)\}, \
0 \leq k \leq K,  \ K \in \Bbb N,
$$
such that
$$
((a^{(k-1)}[-],c^{(k-1)},a^{(k-1)}[+]) , (a^{(k)}[-],c^{(k)},a^{(k)}[+]) ) 
\in \Pi_{\langle H \rangle}(X), \ 0<k \leq K,
$$
and for
$$
r^{(k)} \in\rho((a^{(k-1)}[-],c^{(k-1)},a^{(k-1)}[+]) , (a^{(k)}[-],c^{(k)},a^{(k)}[+]) ), \quad 0 < k \leq K,
$$
it holds that
\begin{multline*}
a^{(0)}[-]^{\mu (a^{(0)}[-])}  
( \prod_{ 0< k \leq K}a^{(k)}[+]^{\mu (a^{(k)}[+])} a^{(k)}[+]^{\gamma(a^{(k)}[+])}_{[1, r^{(k)}]}c^{(k+1)})\\
 a^{(K)}[+]^{\mu (a^{(K}[+])} \in \mathcal L(X),
 \end{multline*}
implies for 
\begin{align*}
\widetilde r^{(k)} \in \Xi (\Theta_{\langle H \rangle}(X,\widetilde X),
(a^{(k-1)}[-],c^{(k-1)},a^{(k-1)}[+]) , (a^{(k)}[-],c^{(k)},a^{(k)}&[+]), r^{(k)}), \\
& 0 < k \leq K,
\end{align*}
 that
\begin{align*}
\widetilde a^{(0)}[-]^{\mu (\widetilde a^{(0)}[-]) }  
 ( \prod_{ 0< k \leq K}\widetilde a^{(k)}[+]^{\mu (\widetilde a^{(k)}[+])} 
 \widetilde a^{(k)}&[+]^{\gamma(\widetilde a^{(k)}[+])}_{[1, \widetilde r^{(k)}]}\widetilde c^{(k+1)})\\
& \widetilde a^{(K)}[+]^{\mu (\widetilde a^{(K}[+])} \in \mathcal L(\widetilde X).
 \end{align*}
 
\begin{lemma}
Let $X$ and $\widetilde X$ be a topologically transitive sofic systems and let there be given a homomorphism $\varphi : X \to \widetilde X $. Let $H \in \Bbb N$ be such that 
$P_{\langle H \rangle}(X) \neq \emptyset$. Then the homomorphism  
$\varphi_\circ =   \varphi \restriction P_{\langle H \rangle}(X)$ has an accompanying map.
\end{lemma}
\begin{proof}
We describe an accompanying map 
$
\Theta_{\langle H \rangle}(X,\widetilde X)
$
for $\varphi_\circ $ by indicating for a given triple 
$$
(a[-],c,a[+]) \in \Omega_{\langle H \rangle}(X) \setminus \{(a, \epsilon, a):
 a  \in \mathcal A_{\langle H \rangle}(X)\}
$$
the structure of its image $(\widetilde a[-],\widetilde c,\widetilde a[+]) $ under 
$\Theta_{\langle H \rangle}(X,\widetilde X)$.
In the case, that 
$$
\varphi(z^{(a[-],c,a[+]) })
\not \in P(\widetilde X),
$$
we apply Lemma (4.1) to obtain
$
 \widetilde t^\prime \in \Bbb Z$ and $\widetilde a^\prime[-],
 \widetilde a[+] \in \mathcal A(\widetilde X)$, such that 
 $$
 \ell (\widetilde a^\prime[-])   =  \ell (^\varphi a[-]),  \qquad   \ell (\widetilde a[+])=
 \ell (^\varphi a[+]),
 $$
 and $\widetilde c^\prime \in \mathcal L(\widetilde X )$,
 such that $\widetilde  z = \varphi(z^{(a[-],c,a[+]) }) $ is given by
 $$
\widetilde z_{(-\infty, \widetilde t^\prime ]}=p^{(-,\widetilde a[-])}\widetilde c^\prime, \quad \widetilde z_{[\widetilde t^\prime, \infty)}=p^{(+,\widetilde a[+])}.
 $$

In the case that $\widetilde c^\prime\neq \epsilon$,
to obtain the complete triple $(\widetilde a[-],\widetilde c,\widetilde a[+]) $, set
$$
\widetilde c =( \Psi_{\widetilde a^\prime[-]} (\widetilde c^\prime))
_{(\kappa_{\widetilde a^\prime[-]}( \Psi_{\widetilde a^\prime[-]}(\widetilde c^\prime)), \ell(\widetilde c^\prime)]},
$$
$$
\widetilde a[-] = (\widetilde a^\prime[-]\Psi_{\widetilde a^\prime[-]} (\widetilde c^\prime))_
{(\kappa_{\widetilde a^\prime[-]}( \Psi_{\widetilde a^\prime[-]}(\widetilde c^\prime)), 
\ell(\widetilde a^\prime[-]) + 
\kappa_{\widetilde a^\prime[-]}( \Psi_{\widetilde a^\prime[-]}(\widetilde c^\prime))]},
$$
 and obtain $\widetilde t$, such that 
 $0 \leq \widetilde t < \beta(\widetilde a[-],\widetilde c,\widetilde a[+])$, 
 by adjusting $\widetilde t^\prime$.
 
 In the case that $\widetilde c^\prime = \epsilon$, set
 $$
 \widetilde c = \epsilon, \qquad \widetilde a[-]= \widetilde a^\prime[-],
 $$
and obtain $\widetilde t$, such that 
 $0 \leq \widetilde t < \beta(\widetilde a[-],\widetilde c,\widetilde a[+])$, 
 by adjusting $\widetilde t^\prime$.
 
 In the case that 
$$
\varphi(z^{(a[-],c,a[+]) })
\in P(\widetilde X),
$$
set
$$
\widetilde t = 0 \quad  \widetilde c = \epsilon, \quad \widetilde a[-] = \widetilde a[+] 
= ^\varphi a[+].
$$ 

Conditions $(A)$ and $(B)$ are satisfied by construction. We indicate how Condition $(B)$ arises.
For this let $[-L, L]$ be a coding window of $\varphi$. Also let $K \in \Bbb N$, and let
$$
(a^{(k)}[-],c^{(k)},a^{(k)}[+]) \in
\Omega_{\langle H \rangle}(X) \setminus \{(a, \epsilon, a): a  \in \mathcal A_{\langle H \rangle}(X)\}, \
0 \leq k \leq K,  \ K \in \Bbb N,
$$
be such that
$$
((a^{(k-1)}[-],c^{(k-1)},a^{(k-1)}[+]) , (a^{(k)}[-],c^{(k)},a^{(k)}[+]) ) \in \Pi_{\langle H \rangle}(X), \ 0<k \leq K,
$$
and also, let
$$
r^{(k)} \in \rho((a^{(k-1)}[-],c^{(k-1)},a^{(k-1)}[+]) , (a^{(k)}[-],c^{(k)},a^{(k)}[+]) ), \quad 0 < k \leq K.
$$
Let
$$
l_k \geq \max_{a\in \mathcal A_{\langle H \rangle}(X)}\ell(a^\varphi )\mu(a^\varphi ), \qquad 0 < k \leq K,
$$ 
and consider the point $x\in X$, that is given by
$$
x_{(-\infty, 0]}=
p^{(-,a^{(0)}[-])}c^{(0)} ,
$$
and
$$
x_{[1, \infty)}=  \left( \prod_{ 0< k < K}a^{(k)}[+]^{\mu (a^{(k)}[+])
+ 2L + l_k} a^{(k)}[+]^{\gamma(a^{(k)}[+])}_{[1, r^{(k)}]}c^{(k+1)}\right)
p^{(+, a^{(K)}[+])}.
$$
A comparison of the image of $x$ under $\varphi$ with the images of the points 
$$
z^{(a^{(k)}[-],c^{(k)},a^{(k)}[+])}, \qquad 0 < k \leq K,
$$ 
under $\varphi$
shows, that $(B)$ holds.
\end{proof}

As a consequence of the topological transitivity of $X$ one has that
$$
\lim _{l \to \infty} \card (( \mathcal L_l(X) \cap \delta)) = \infty, \qquad \delta \in \Delta.
$$
Denote by $H_0(X)$ the smallest $H_0\in\Bbb N$ such that 
$$
\card (\{u \in  \mathcal L_{H_0}(X): [u]_\approx \in \delta\}    ) > \card(\Delta(X)) \max\{ \ell(f( \delta ):  
\delta \in  (\Delta( X)\}.
$$

Given a topologically transitive sofic system $X$ we use the notation
$$
\lambda(X) = \max \{\eta (\delta): \delta \in \Delta \}(\max \{\nu (\delta): \delta \in \Delta \} + \max \{\eta (\delta): \delta \in \Delta \}),
 $$
 $$
 H_1(X) =  \lambda(X) + \max \{\ell(f(\delta))\mu (f(\delta))\gamma(f(\delta)): \delta \in \Delta \}.
$$
Given a topologically transitive sofic systems $X$ and $\widetilde X $ we use the notation
$$
H_2(X, \widetilde X) =   \lambda(\widetilde X)  \max_{ \delta \in \Delta} \{  \max_{\{\widetilde a \in \mathcal A_{\langle H \rangle}(\widetilde X) : \ell(\widetilde a) \vert  \ell(a) \}}
   \ell(\widetilde a) \mu(\widetilde a)\gamma(\widetilde a)   \}.
$$
\begin{lemma} 
Let $X$ and $\widetilde X$ be  topologically transitive sofic systems. Let
\begin{align*}
H > H_0(X) + \max \{H_1(X), H_2(X) \}, \tag {4.2}
\end{align*}
and let there be given a homomorphism
$$
\varphi_\circ : P_{\langle H \rangle}(X) \to P_{\langle H \rangle}(\widetilde X)
$$
with an accompanying map $\Theta_{\langle H \rangle}(X,\widetilde X)$. Then there exists a homomorphism 
$$
\varphi : X \to \widetilde X
$$
that extends $\varphi_\circ.$
\end{lemma}
\begin{proof} 
Let $F\subset X$ be such that (2.2) and (2.3) hold. 
Set
\begin{align*}
\widehat H = \max_{a \in \mathcal A_{\langle H \rangle}(X)} \ell(a) \mu(a)\gamma(a) +
 \max_{a \in \mathcal A_{\langle H \rangle}(X)}
  \{  \max_{\{\widetilde a \in \mathcal A_{\langle H \rangle}(\widetilde X) : \ell(\widetilde a) \vert  \ell(a) \}}
   \ell(\widetilde a) \mu(\widetilde a)\gamma(\widetilde a)   \}. 
   \tag{4.3}
\end{align*}
and
\begin{align*}
\widehat F = F \cup S_X^H \left(F_\infty^+  \cup   \bigcup_{3H \leq j < \infty }  F_j^+\right)
\ \cup \  &S_X^{-H}\left(F_\infty^- \cup \bigcup_{3H \leq j < \infty }  F_j^-\right) 
\tag {4.4}
\  \cup  \
\\
& \left (\bigcup_{4H \leq j < \widehat H} (\bigcup_{1 < k <  \frac {j}{H}} S_X^{kH}(F_j^+ ))\right).
\end{align*}
We set
$$
\widehat F^- = \bigcup_{H \leq J < \infty} \widehat F^-_J, \quad 
\widehat F^+ = \bigcup_{H \leq J < \infty} \widehat F^+_J.
$$
One has that
\begin{align*}
\widehat F \setminus \widehat F^- \subset
 \left(\bigcup_{ \widehat H \leq J < \infty}  \widehat F^-_J  \right) \cup \widehat F^-_\infty, 
 \quad
 \widehat F \setminus \widehat F^- \subset
 \left(\bigcup_{ \widehat H \leq J < \infty}  \widehat F^-_J  \right) \cup \widehat F^-_\infty.
  \tag{4.5}
\end{align*}

For $ \delta \in \Delta$ we choose a word
$
b(\delta) \in 
\delta \cap 
\mathcal L_{H_0}(X)
$
such that
\begin{align*}
b(\delta)\in \delta \setminus  
\{p_{[1 + l, H_0 + l]}
^{(+, f(\delta))}): 0 < l < \ell (f(\delta))\},    
\end{align*}
and for $u \in\mathcal L_{H_0}(X)$ we  choose a word
$
b(u) \in  \mathcal L_{H_0}(X)
$
such that
\begin{align*}
 b(u) \in [u]_\approx\setminus\{u\}.
\end{align*}

By shift invariance we obtain an endomorphism $x \to \bar {x}\in X  \ (x \in X)$ by replacing simultaneously
in the points $x \in \widehat F^+ $,  
the block $x_{[1, H_0]} $ by a block that carries the word $b( x_{[1, H_0]} )$,
 and in the points $x \in \widehat F \setminus  \widehat F^+ $
 the block $x_{[H_0, \widehat i^+(x)]}$ 
 by a block, that carries the word 
 $b(\delta(x_{[H_0, \widehat i^+(x)]}) )$.
 
 We set
 $$
 I^-(x) = I^-(x_{(H_0, \widehat i^+(x)]}),  \  f(x) = f(x_{(H_0, \widehat i^+(x)]}), \qquad
 x \in \widehat F^+.
 $$
In view of (4.3)  and (4.4)  one has
$$
x_{( -2H, 0]} \in \mathcal P_{2H}, \qquad x\in \widehat F \setminus \widehat F^-,
$$
$$
x_{[1, 2H]} \in \mathcal P_{2H}, \qquad x\in \widehat F \setminus  \widehat F^+,
$$
and we set
$$
g_-(x) = g (x_{( -2H, 0]} )  \qquad x\in \widehat F \setminus \widehat F^-,
$$
$$
g_+(x) = g (x_{[1, 2H]} )  \qquad x\in \widehat F \setminus \widehat F^+. 
$$
By 
the choice of the words 
$b(u), u \in\mathcal L_{H_0}(X),$
 we can at this stage  secure for $x\in \widehat F $  a time point $T(x) \in \Bbb N$  by
 \begin{align*}
T(x) =
\begin{cases}
\min\{i \in [1, H_0 + I^-(x)]:  \bar x_{(i , H_0  + I^-(x), 0]} =  \tag{4.6}
 p^{(-, f( x))}_{( i - H_0  - I^-(x), 0]}\}, &\text{if  $x \in \widehat F^+ $}, \\
\min\{i \in [1, H_0 ]:  \bar x_{(i , H_0]} = 
 p^{(-, g_+( x))}_{( i - H_0  , 0]}\}, &\text {if $x\in \widehat F \setminus \widehat F^+$}.
\end{cases}  
\end{align*}
We set
\begin{align*}
T^+(x) =T(S^{\widehat i^+(x)}_X(x)),  
\
T^-(x) =T(S^{-\widehat i^-(x)}_X(x)), \qquad x \in \widehat F.
\end{align*}
One has
\begin{align*}
I^-(x) = I^-(x_{H_0,T^+(x)}),  \  f(x) = f(x_{H_0,T^+(x)} ), \qquad  x \in \widehat F^+. 
\tag {4.7}
\end{align*}
We  set
$$
a(x) = (x_{[1, H_0]}\Psi(x_{[ H_0, T^+(x) ]})_{(T(x), T(x) +   \ell(f(x)]}, \qquad 
x \in \widehat F^+.
$$
By the choice of $H_1(X)$ and by Lemma 3.1 
\begin{align*}
 (x_{[1, H_0]}\Psi(x_{[ H_0, T^+(x) ]})_{(T(x), T^+(x) ]} 
 \in \Gamma^+( a(x)^{\mu(a(x))\gamma (a(x) )} ), \qquad   x\in\widehat F^+, \tag {4.8}
\end{align*}
and by (4.3) and (4.4)
\begin{align*}
x_{[1, T(x)]}\in  \Gamma^+(g_-(x)^{\mu( g_-(x))}), \qquad 
x\in \widehat F \setminus \widehat F^+. \tag {4.9}
\end{align*}
With the use of the auxiliary maps $\Psi$ and $\Psi_a$ we define by shift invariance an endomorphism
 $x \to \bar {\bar x} \in  X \  (x \in X)$ by replacing simultaneously for 
  $x\in \widehat F^+$  in the points $\bar x$ the blocks $\bar x_{(T(x), T^+(x)]}$ by blocks that carry the word
  $$
  \Psi_{a(x)} (x_{[1, H_0)]}(\Psi(x_{(H_0, T^+(x) ]})_{(T(x), T^+(x) ]} ),
  $$
  and for $x\in \widehat F \setminus \widehat F^+$ by replacing  the blocks $x_{[1, T(x)]}$ by blocks that carry the word
  $$
  \Psi_{g_-(x)}(x_{[1, T(x)]}).
  $$
 By (4.8) the first part of this coding instruction is meaningful, and by (4.9) its second part is also meaningful. By (4.7) its first part is also consistent with the definition of $T^+$.
 
 At this stage triples
\begin{align*}
((a[-],c,a[+]) )(x) \in  \Omega_{\langle H \rangle}(X)
 \setminus \{(a, \epsilon, a): a  \in \mathcal Q_{\langle H \rangle}(X)\}), \qquad x\in \widehat F,
 \tag {4.10}
\end{align*}
 have become visible.
We set
$$
a[+](x) =
\begin{cases}
a(x), &\text{if  $x\in \widehat F$}, \\
g(x_{(T(x), T(x) + 2H ]}), &\text {if $x\in \widehat F \setminus \widehat F^+$},
\end{cases}
$$
and with the notation
$$
l(x) =
\begin{cases}
T(x) - T^-(x) -\kappa_{g(x)}(  \bar {\bar x}_{(T^-(x), T(x)]}  ), &\text{if  $x\in \widehat F$}, \\
T(x) -  \kappa_{a( S_X^{\widehat i^+(x)})}(  \bar {\bar x}_{[1, T(x)]}) , &\text {if $x\in \widehat F \setminus \widehat F^+$}.
\end{cases}
$$
$$
l^{\prime}(x) =
\begin{cases}
\ell(a( S_X^{\widehat i^+(x)})), &\text{if  $x\in \widehat F$}, \\
\ell (g(x)), &\text {if $x\in \widehat F \setminus \widehat F^+$},
\end{cases}
$$
we set
$$
c(x) =  \bar {\bar x}_{[T^-(x)- l(x),T(x)]},
\
a[-](x) =  \bar {\bar x}_{[T^-(x)- l(x) - l^{\prime}(x) ,T(x) - l(x]}, \ \  x \in \widehat F.
$$
The maps $\Psi_a$ do not change the final symbol of a word. Also by its definition $\mu$ only depends on  conjugacy classes, and it is seen that (4.10) holds.

With the triples $((a[-],c,a[+]) )(x), x \in \widehat F$, at hand the structure of the words that are carried by the blocks 
$ \bar {\bar x}_{(T(x),T^+(x)]}, x \in \widehat F$,
has become apparent. With the notation
$$
q(x) = \lfloor T^+(x) - T(x) - \ell(c(S_X^{\widehat i^+(x)}(x))/  
 \ell (a[+](x))\gamma(a[+](x))\rfloor,
$$
$$
r(x) = T^+(x) - T(x) - \ell(c(S_X^{\widehat i^+(x)}(x)) -  q(x)\ell (a[+](x))\gamma(a[+](x)),
$$
these words are given by
\begin{align*}
 \bar {\bar x}_{(T(x),T^+(x)]} = 
 a[+](x)^{q(x)\gamma(a[+](x))}a[+](x)^{\gamma(a[+](x))}_{[1, r(x)]}
 c(S_X^{\widehat i^+(x)}(x)). \tag{4.11}
 \end{align*}
 Also
 \begin{align*}
  \bar {\bar x}_{(-\ell(c(x),0]} = c(x), \qquad x \in  \widehat F^-. \tag{4.12}
 \end{align*}
 Also note that the block that the blocks $\bar {\bar x}_{(-\infty,-\ell(c(x)))]}, 
 x \in \widehat F^-,$ carry the semi-infinite concatenation of $a[-](x)$, and that the blocks 
 $\bar {\bar x}_{[1,\infty)}, 
 x \in \widehat F^+,$ carry the semi-infinite concatenation of $a[+](x)$.

By mans of the accompanying map $\Theta_{\langle H \rangle}(X,\widetilde X)$ we define a homomorphism
$
\varphi: X \to  \widetilde X.
$
We denote the image of the triple $(a[-],c,a[+])(x), x \in \widehat F$, under $\Theta_{\langle H \rangle}(X,\widetilde X)$ by
$
((\widetilde a[-],\widetilde c,\widetilde a[+]) ,\widetilde t(x)),
$
 and we set
$$
\widetilde  T(x) = T(x) + \widetilde t(x), \  \widetilde T^+(x) = T^+(x) +
\widetilde t(S_X^{i^+_{\widehat F}}(x)), \quad
x \in  \widehat F.
$$
Following the pattern set by (4.11) and (4.12) we use the notation
$$
 \widetilde q(x) = \lfloor  \widetilde T^+(x) - \widetilde T(x) - 
 \ell( \widetilde c(S_X^{i^+_{\widehat F}}(x))/  
\gamma( \widetilde a[+](x)) \ell ( \widetilde a[+](x))\rfloor,
$$
\begin{multline*}
 \widetilde r(x) =  \widetilde T^+(x) -  \widetilde T(x) - 
 \ell( \widetilde c(S_X^{ i_{\widehat F}^+(x)}(x)) - \widetilde  q(x) \gamma( \widetilde a[+](x))\ell ( \widetilde a[+](x)),  
 \\
  x \in  \widehat F\setminus  \widehat F^+.
\end{multline*}
and we set
\begin{multline*}
 \widetilde x_{[\widetilde T(x),\widetilde T^+(x)]} = 
 \widetilde a[+](x)^{\widetilde q(x)\gamma(\widetilde a[+](x))}\widetilde a[+](x)^{\gamma
 \widetilde a[+](x))}_{[1, \widetilde r(x)]}
 \widetilde c(S_X^{i^+_{\widehat F}}(x)), 
 \\
  x \in  \widehat F\setminus  \widehat F^+.
 \end{multline*}
 
We set
$$
\widetilde x_{(\widetilde T(x) - \ell(\widetilde c(x)),\widetilde T(x)]} = \widetilde c(x), \qquad
x \in  \widehat F^-.
$$
We denote the block map with coding window $[-H, H]$, that implements the homomorphism $\varphi_\circ$ by $\bold {\Phi_\circ}$, and  for $x\in X$, such that
$$
i^+_{\widehat F}(x) \geq  \ell(\widetilde c(x)) - \widetilde T(S_X^{i^+_{\widehat F}(x)}(x)),
$$
and
$$
i^-_{\widehat F}(x) \geq \widetilde T(S_X^{-i^-_{\widehat F}(x)}(x)),
$$
we set
$$
\widetilde x_0 = \bold \Phi_\circ (x_{[- H, H]}).   \qed
$$
\renewcommand{\qedsymbol}{}
\end{proof}

\section{Decidability}

 Lemma (4.2) and Lemma (4.3) together contain a necessary and sufficient condition on topologically transitive sofic systems  $X$ and $\widetilde X$ for the existence of a homomorphism $\varphi:X\to\widetilde X$. This condition is decidable, as the mathematical structures that enter  are finite (see  Lemma (3.1)). Also there is a bound on the parameter $K$ of Condition $(B)$. To see this, assume that there is no homomorphism from $X$ to $\widetilde X$, and denote
by $K_{min}$ the minimal $K\in \Bbb N$ such that there are
$$
(a^{(k)}[-],c^{(k)},a^{(k)}[+]) \in
\Omega_{\langle H \rangle}(X) \setminus \{(a, \epsilon, a): a  \in \mathcal A_{\langle H \rangle}(X)\}, \
0 \leq k \leq K_{min},  \ K \in \Bbb N,
$$
such that
$$
((a^{(k-1)}[-],c^{(k-1)},a^{(k-1)}[+]) , (a^{(k)}[-],c^{(k)},a^{(k)}[+]) ) 
\in \Pi_{\langle H \rangle}(X), \ 0<k \leq K_{min},
$$
and 
$$
r^{(k)} \in((a^{(k-1)}[-],c^{(k-1)},a^{(k-1)}[+]) , (a^{(k)}[-],c^{(k)},a^{(k)}[+]) ), \quad 0 < k \leq K_{min},
$$
such that 
\begin{multline*}
a^{(0)}[-]^{\mu (a^{(0)}[-])}  
( \prod_{ 0< k \leq K}a^{(k)}[+]^{\mu (a^{(k)}[+])} 
a^{(k)}[+]^{\gamma(a^{(k)}[+])}_{[1, r^{(k)}]}c^{(k+1)})\\
 a^{(K_{min})}[+]^{\mu (a^{(K_{min}}[+])} \in \mathcal L(X),
 \end{multline*} 
and
\begin{align*}
\widetilde r^{(k)} \in \Xi (\Theta_{\langle H \rangle}(X,\widetilde X) ,((a^{(k-1)}[-],c^{(k-1)},a^{(k-1)}[+]) , (a^{(k)}[-],c^{(k)},a^{(k)}&[+]) ), r^{(k)}), \\
& 0 < k \leq K_{min},
\end{align*}
such that
\begin{align*}
\widetilde a^{(0)}[-]^{\mu (\widetilde a^{(0)}[-]) }  
 ( \prod_{ 0< k < K}\widetilde a^{(k)}[+]^{\mu (\widetilde a^{(k)}[+])} 
 \widetilde a^{(k)}&[+]^{\gamma(\widetilde a^{(k)}[+])}_{[1, \widetilde r^{(k)}]}\widetilde c^{(k+1)})\\
& \widetilde a^{(K_{min})}[+]^{\mu (\widetilde a^{(K_{min}}[+])} \not\in \mathcal L(\widetilde X). \tag {5.1}
 \end{align*}
With $H\in \Bbb N$ as in Lemma (4.2) one has that
$$
K_{min} \leq \card ( \Delta(X)) \card (  \Delta(\widetilde X) ) \card (\Omega (X)) 
\card (\Omega (\widetilde X)) H^2.
$$
Otherwise, there would exist $\bar K,  \bar{\bar K} \in [1,  K_{min}]$ such that
\begin{multline*}
(a^{(\bar K)}[+], r^{(\bar K)},c^{(\bar K+1)},\widetilde a^{(\bar K)}[+], \widetilde r^{(\bar K)},\widetilde c^{(\bar K+1)} )= \\
(a^{(\bar{\bar K})}[+], r^{(\bar{\bar K})},c^{(\bar{\bar K}+1)},\widetilde a^{(\bar{\bar K})}[+], \widetilde r^{(\bar{\bar K})},\widetilde c^{(\bar{\bar K}+1)} ),\tag{5.2}
\end{multline*}
\begin{align*}
&[a^{(0)}[-]^{\mu (a^{(0)}[-])}  
( \prod_{ 0< k \leq \bar K}a^{(k)}[+]^{\mu (a^{(k)}[+])} a^{(k)}[+]^{\gamma(a^{(k)}[+])}_{[1, r^{(k)}]}c^{(k+1)})]_\approx =\tag {5.3}
 \\
&[a^{(0)}[-]^{\mu (a^{(0)}[-])}  
( \prod_{ 0< k \leq \bar{\bar K} }a^{(k)}[+]^{\mu (a^{(k)}[+])} a^{(k)}[+]^{\gamma(a^{(k)}[+])}_{[1, r^{(k)}]}c^{(k+1)})]_\approx,
\end{align*}
\begin{align*}
&[\widetilde a^{(0)}[-]^{\mu (\widetilde a^{(0)}[-])}  
( \prod_{ 0< k \leq \bar K}\widetilde a^{(k)}[+]^{\mu (\widetilde a^{(k)}[+])} \widetilde a^{(k)}[+]^{\gamma(\widetilde a^{(k)}[+])}_{[1,\widetilde r^{(k)}]}c^{(k+1)})]_\approx = \tag{5.4}
\\
&
[\widetilde a^{(0)}[-]^{\mu (\widetilde a^{(0)}[-])}  
( \prod_{ 0< k \leq \bar{\bar K}}\widetilde a^{(k)}[+]^{\mu (\widetilde a^{(k)}[+])} 
\widetilde a^{(k)}[+]^{\gamma(\widetilde a^{(k)}[+])}_{[1,\widetilde r^{(k)}]}c^{(k+1)})]_\approx.
\end{align*}
Set
$$
K_{\circ} = K_{min}- \bar{\bar K}+{\bar K}.
$$
By the minimality  assumption on $K_{min}$, and by  (5.2), (5.3), and (5.4)  one obtains by
collapsing the sequence 
$$
(a^{(k)}[+]), r^{(k)},c^{(k+1)},\widetilde a^{(k)}[+]),\widetilde  r^{(k)},\widetilde c^{(k+1)}  )_
{1 \leq k <  K_{min}},
$$
by
\begin{multline*}
(a_{\circ}^{(k)}[+]), r_{\circ}^{(k)},c_{\circ}^{(k+1)} ,\widetilde a_{\circ}^{(k)}[+],
\widetilde r_{\circ}^{(k)},\widetilde c_{\circ}^{(k+1)} ) = \\
(a^{(k)}[+]), r^{(k)},c^{(k+1)},\widetilde a^{(k)}[+]),\widetilde r^{(k)},\widetilde c^{(k+1)}) ,
\qquad 1 \leq k \leq \bar K,
\end{multline*}
and
\begin{multline*}
(a_{\circ}^{(k)}[+]), r_{\circ}^{(k)},c_{\circ}^{(k+1)} ,\widetilde a_{\circ}^{(k)}[+],
\widetilde r_{\circ}^{(k)},\widetilde c_{\circ}^{(k+1)} ) = \\
(a^{(k +\bar{\bar K}-{\bar K})}[+], r^{(k +\bar{\bar K}-{\bar K})},c^{(k +\bar{\bar K}-{\bar K}+1)},\widetilde a^{(k +\bar{\bar K}-{\bar K})}[+],\widetilde r^{(k +\bar{\bar K}-{\bar K})},\widetilde c^{(k +\bar{\bar K}-{\bar K}+1)}  ), \\ 
\bar K < k < K_{\circ},
\end{multline*}
a sequence
$(a_{\circ}^{(k)}[+], r_{\circ}^{(k)},c_{\circ}^{(k+1)},\widetilde a_{\circ}^{(k)}[+],\widetilde r_{\circ}^{(k)},\widetilde c_{\circ}^{(k+1)})_{1 \leq k <  K_{\circ}},$
such that

\begin{align*}
 a^{(0)}[-]^{\mu ( a^{(0)}[-]) }  
 \bigr( \prod_{ 0< k < K_{\circ}} a_{\circ}^{(k)}[+]^{\mu ( a_{\circ}^{(k)}[+])} 
 a_{\circ}^{(k)} &[+]^{\gamma(a_{\circ}^{(k)}[+])}_{[1, r_{\circ}^{(k)}]} c_{\circ}^{(k+1)}
  \bigr)
 \\
& a_{\circ}^{(K)}[+]^{\mu ( a_{\circ}^{(K_{\circ})}[+])}\in \mathcal L( X),
 \end{align*}
\begin{align*}
\widetilde a^{(0)}[-]^{\mu (\widetilde a^{(0)}[-]) }  
( \prod_{ 0< k < K_{\circ}}\widetilde a_{\circ}^{(k)}[+]^{\mu (\widetilde a_{\circ}^{(k)}[+]) 
} \widetilde a_{\circ}^{(k)}&[+]^{\gamma(\widetilde a_{\circ}^{(k)}[+])}_{[1,\widetilde r_{\circ}^{(k)}]}
\widetilde c_{\circ}^{(k+1)})
 \\
& \widetilde a_{\circ}^{(K)}[+]^{\mu (\widetilde a_{\circ}^{(K)}[+])}\in \mathcal L(\widetilde X),
 \end{align*}
which then by (5.4) contradicts (5.1).

\section{Surjective homomorphisms}
 We consider a topologically mixing sofic system $X$. For $H\in \Bbb N$ we denote by 
 $\Omega^-_{\langle  H \rangle}(X)$
the set of pairs
$$
(\delta, a) \in \Delta_{synchro}(X) \times \mathcal A_{\langle  H \rangle} (X)
$$
such that 
$$
\delta  a^{\mu (a)}\subset \mathcal L(X),
$$
and such  that $\delta$ contains a word with a final symbol, that is different from the final symbol of $a$, and we denote by 
 $\Omega^+_{\langle  H \rangle}(X)$
the set of pairs
$$
(a, \delta) \in  \mathcal A_{\langle  H \rangle} (X)\times\Delta_{synchro}(X)
$$
such that 
$$
a^{\mu (a) }\delta  \subset \mathcal L(X),
$$
and such that $\delta$ contains a word with an initial symbol, that is different from the initial symbol of $a$.

For a pair
$$
((\delta, a), (a[-],c,a[+])) \in \Omega^-_{\langle  H \rangle}(X)\times 
(\Omega_{\langle H \rangle}(X) \setminus \{(a, \epsilon, a): a  \in \mathcal A_{\langle H \rangle}(X)\}),
$$
such that $a$ is conjugate to $a[-]$, we denote by
$$
\rho^-((\delta, a), (a[-],c,a[+]))
$$
the set of $r\in [1, \ell(a)]$, such that the initial symbol of $\tau^{(r)}(a)$ is different from the initial symbol of $c$ (in other words, such that 
$((\delta, a),(\tau^{(r)}(a),c,a[+])) \in  \Omega^-_{\langle  H \rangle}(X)$). The meaning of $\rho^+$ is symmetric. We denote by
$\Pi^-_{\langle  H \rangle}(X)$ the set of pairs 
$$
((\delta, a), (a[-],c,a[+])) \in \Omega^-_{\langle  H \rangle}(X)\times 
(\Omega_{\langle H \rangle}(X) \setminus \{(a, \epsilon, a): a  \in \mathcal A_{\langle H \rangle}(X)\}),
$$
such that $a$ is conjugate to $a[-]$, and such that
$$
a[-] = \tau^{(\min \rho^-((\delta, a), (a[-],c,a[+])))}(a).
$$
The meaning of $\Pi^+_{\langle  H \rangle}(X)$ is symmetric.

Assume also given a sofic system $\widetilde X$. We consider a map
$$
\Theta^-_{\langle H \rangle}(X,\widetilde X):  \Omega^-_{\langle  H \rangle}(X) \to
 \Omega^-_{\langle  H \rangle}(\widetilde  X) \times \Bbb Z_+
$$
that assigns to a pair 
$
(\delta, a) \in  \Omega^-_{\langle  H \rangle}(X) 
$
an element
$
((\widetilde a , \widetilde\delta), \widetilde t^-)$ of $ \Omega^-_{\langle  H \rangle}(\widetilde  X) \times \Bbb Z_+
$
such that
\begin{align*}
0\leq  \widetilde t^- <\gamma(\widetilde a) \ell(\widetilde a )). 
\end{align*}

Given  
$\Theta^-_{\langle H \rangle}(X,\widetilde X)$  denote for
$$
((\delta(-), a(-)),(a[-],c,a[+])  ) \in \Pi^-_{\langle H \rangle}(X)
$$
the image under $\Theta^-_{\langle H \rangle}(X,\widetilde X)$  of 
$(\delta(-), a(-))$ by  $( \widetilde\delta(-), \widetilde a(-)), \widetilde t^-)$
and the image under $\Theta^-_{\langle H \rangle}(X,\widetilde X)$ of
$(a[-],c,a[+])$ by
$((\widetilde a[-],\widetilde c,\widetilde a[+]), \widetilde t)$, 
and we denote for
$$
r^- \in  \rho^-((\delta(-), a(-)), (a[-],c,a[+]))
$$
by 
$ \Xi^-(\Theta^-_{\langle H \rangle}(X,\widetilde X),(\delta(-), a(-)), (a[-],c,a[+]), r^-)$ the set of remainders that are left by the numbers in
$$
\{r + \ell(c) - \ell(\widetilde c) + \widetilde t^- - \widetilde t + m: 
0\leq m <\gamma(\widetilde a) \ell(\widetilde a ) \}
$$
after division by
$\gamma(\widetilde a) \ell(\widetilde a).$ The meaning  $\Xi^+$ of is symmetric.
Given also $H\in \Bbb N$, such that $P_{\langle H \rangle}\neq \emptyset$
and a homomorphism
$$
\varphi_\circ: P_{\langle  H \rangle}(X) \to P_{\langle  H \rangle}(\widetilde X),
$$
we say that the triple 
$
(\Theta^-_{\langle H \rangle}, \Theta_{\langle H \rangle},
\Theta^+_{\langle H \rangle})(X,\widetilde X)
$
accompanies $\varphi_\circ $, if $\Theta_{\langle H \rangle}(X,\widetilde X) $ accompanies $\varphi_\circ $, and if the following statements $A(-) $,  
$B(-) $ and $ B(0) $
hold, as well as a statement $A(+) $, that is symmetric to $A(-) $, and a statement $ B(+) $, that is symmetric to  $B(-) $:

\noindent
$A(-)$: For $(\delta(-), a(-)) \in \Omega^-_{\langle  H \rangle}(X)$  one has that 
$\widetilde a(-) = \tau^{(\widetilde t^-)}(^\varphi a)$.

\noindent
$B(-)$: For
$$
(\delta(-), a(-)) \in \Omega^-_{\langle  H \rangle}(X),
$$
$$
(a^{(k)}[-],c^{(k)},a^{(k)}[+]) \in  
\Omega_{\langle H \rangle}(X) \setminus \{(a, \epsilon, a): a  \in \mathcal A_{\langle H \rangle}(X)\}, \
0 <k < K,  \ K \in \Bbb N,
$$
such that
$$
((\delta(-), a(-)),(a^{(1)}[-],c^{(1)},a^{(1)}[+])   )  \in \Pi^-_{\langle H \rangle}(X),
$$
$$
((a^{(k-1)}[-],c^{(k-1)},a^{(k-1)}[+]) , (a^{(k)}[-],c^{(k)},a^{(k)}[+]) ) 
\in \Pi_{\langle H \rangle}(X), \ 0<k \leq K,
$$
and for 
$$
r^- \in \rho^-((\delta^-, (a)), (a^{(1)}[-],c^{(1)},a^{(1)}[+])), 
$$
$$
r^{(k)} \in \rho((a^{(k-1)}[-],c^{(k-1)},a^{(k-1)}[+]), (a^{(k)}[-],c^{(k)},a^{(k)}[+])),
 \quad 0 < k \leq K,
$$
it holds that
\begin{multline*}
\delta a^{\mu(a)}a^{\gamma(a)}_{[1, r^-]} c^{(1)} 
\left ( \prod_{ 0< k < K}a^{(k)}[+]^{\mu (a^{(k)}[+])} 
a^{(k)}[+]^{\gamma(a^{(k)}[+])}_{[1, r^{(k)}]}c^{(k+1)}\right)
\\
 a^{(K)}[+]^{\mu (a^{(K)}[+])} \subset \mathcal L(X),
 \end{multline*}
implies for 
$$
\widetilde r^- \in  \Xi^-(\Theta^-_{\langle H \rangle}(X,\widetilde X),(\delta^	-, a(-)), (a[-],c,a[+])), r^-),
$$
\begin{align*}
\widetilde r^{(k)} \in \Xi( \Theta^-_{\langle H \rangle}(X,\widetilde X),(a^{(k-1)}[-],c^{(k-1)},a^{(k-1)}[+]) , (a^{(k)}[-],c^{(k)},a^{(k)}&[+]), r^{(k)}), 
\\
& 0 < k \leq K,
\end{align*}
 that
\begin{multline*}
\widetilde\delta(-)\widetilde a(-)^{\mu(\widetilde a(-))}\widetilde a(-)^{\gamma(\widetilde a(-))}_{[1,\widetilde r^-]} c^{(1)} 
\left( \prod_{ 0< k \leq K}\widetilde a^{(k)}[+]^{\mu (\widetilde a^{(k)}[+]) 
} \widetilde a^{(k)}[+]^{\gamma(\widetilde a^{(k)}[+])}_{[1,\widetilde r^{(k)}]}\widetilde c^{(k+1)}\right)\\
 \widetilde a^{(K)}[+]^{\mu (\widetilde a^{(K)}[+])} \subset \mathcal L(X).
 \end{multline*}
\noindent
$B(0)$: For
$$
(\delta(-), a^{(-)}) \in \Omega^-_{\langle  H \rangle}(X), 
$$
$$
(a^{(k)}[-],c^{(k)},a^{(k)}[+]) \in  
\Omega_{\langle H \rangle}(X) \setminus \{(a, \epsilon, a): a  \in \mathcal A_{\langle H \rangle}(X)\}, \
0 <k \leq K,  \ K \in \Bbb N,
$$
$$
(a^{(+)}, \delta(+)) \in \Omega^+_{\langle  H \rangle}(X),
$$
such that
$$
((\delta(-), a^{(-)}),(a^{(1)}[-],c^{(1)},a^{(1)}[+])   ) \in \Pi^-_{\langle H \rangle}(X), 
$$
$$
((a^{(k-1)}[-],c^{(k-1)},a^{(k-1)}[+]) , (a^{(k)}[-],c^{(k)},a^{(k)}[+]) ) 
\in \Pi_{\langle H \rangle}(X), \ 0<k \leq K,
$$
$$
((a^{(K-1)}[-],c^{(K-1)},a^{(K-1)}[+]),(a(+), \delta(+)))  \in \Pi^+_{\langle H \rangle}(X),
$$
and for 
$$
r^- \in \rho^-((\delta(-), a^{(-)}) , (a^{(1)}[-],c^{(1)},a^{(1)}[+]) ), 
$$
$$
r^{(k)} \in \rho((a^{(k-1)}[-],c^{(k-1)},a^{(k-1)}[+]) , (a^{(k)}[-],c^{(k)},a^{(k)}[+]) ), 
\quad 0 < k \leq K,
$$
$$
r^+ \in \rho^+((a^{(K)}[-],c^{(K)},a^{(K)}[+]),(a(+), \delta(+)), 
$$
it holds that
\begin{multline*}
\delta(-)a(-)^{\mu(a(-))}a(-)^{\gamma(a(-))}_{[1, r^-]} c^{(1)} 
 \left( \prod_{ 0< k < K}a^{(k)}[+]^{\mu (a^{(k)}[+])} a^{(k)}[+]^{\gamma(a^{(k)}[+)}_{[1, r^{(k)}]}c^{(k+1)}\right)
 \\
 a(+)^{\mu (a(+))}a(+)^{\gamma(a(+)}_{[1, r^+]}\delta^+ \subset \mathcal L(X),
 \end{multline*}
implies for 
$$
\widetilde r^- \in \Xi^{-}(\Theta^-_{\langle H \rangle}(X,\widetilde X),(
\delta^-, a(-)), (a[-],c,a[+]), r^-),
$$
\begin{align*}
\widetilde r^{(k)} \in \Xi (\Theta^-_{\langle H \rangle}(X,\widetilde X) ,
(a^{(k-1)}[-],c^{(k-1)},a^{(k-1)}[+]) , (a^{(k)}[-],c^{(k)},a^{(k)}&[+]), r^{(k)}), 
\\
 &0 < k < K,
\end{align*}
$$
\widetilde r^+ \in \Xi^{+}(\Theta^+_{\langle H \rangle}(X,\widetilde X) , (a[-],c,a[+]),
(a(+),\delta(+), r^+),
$$
that
\begin{multline*}
\widetilde\delta^-\widetilde a^{\mu(\widetilde a(-))}\widetilde a(-)^{\gamma(\widetilde a(-)}_{[1, \widetilde r^-]}
\widetilde c^{(1)} 
 \left( \prod_{ 0< k < K}\widetilde a^{(k)}[+]^{\mu (\widetilde a^{(k)}[+])
} \widetilde a^{(k)}[+]^{\gamma( \widetilde a^{(k)}}_{[1, \widetilde r^{(k)}]}\widetilde c^{(k+1)}\right)
 \\
 \widetilde a(+)^{\mu (\widetilde a(+))}\widetilde a(+)_{[1, \widetilde r^+]}\widetilde
 \delta(+) \subset \mathcal L(\widetilde X).
 \end{multline*}
\medskip

For the proof of next lemma compare the proof of Lemma 4.2.

\begin{lemma} 
Let $X\subset \Sigma^{\Bbb Z}$ and $\widetilde X\subset\widetilde\Sigma^{\Bbb Z}$ be a topologically mixing sofic systems and let there be given a surjective homomorphism $\varphi : X \to \widetilde X $. 
Let $H \in \Bbb N$, such that $P_{\langle  H \rangle} (X) \neq \emptyset$. Then the homomorphism 
$\varphi_\circ = \varphi\restriction P_{\langle  H \rangle} (X)$ has an accompanying triple
$( \Theta^-_{\langle H \rangle} , \Theta_{\langle H \rangle} , \Theta^+_{\langle H \rangle} )(X,\widetilde X) $.
\end{lemma}
\begin{proof} 
Let $[-L,+L]   $ be a coding window for  $\varphi$ and let $\varphi$ be given by the block map 
$\bold\Phi: X_{[-L,+L]}  \to \widetilde\Sigma$. The map 
$ \Theta_{\langle H \rangle}(X,\widetilde X)$ is as constructed in Section 4.

The construction of the map $\Theta^-_{\langle H \rangle}(X,\widetilde X)$ necessitates  certain choices, that have to be made beforehand. These choices  are possible under the hypothesis that $\widetilde X$  topologically mixing,  and that $\varphi$ is  surjective. For 
$a \in \mathcal A(X)$ we choose a symbol $\widetilde \sigma_a \in\widetilde \Sigma$, that is different from the last symbol of $^\varphi a$, and we choose a word $c_a \in \mathcal L_{2L + 1}$, such that
 $
\bold \Phi(c_a) = \widetilde \sigma_a .
 $
 For $(\delta, a) \in \Omega^-_{\langle  H \rangle}(X)$ we choose a word 
 $d_{(\delta, a)}\in \Gamma^+(c_a)$, that has a suffix in $\delta$, and such that $\ell(d_{(\delta, a)} ) +2L$ is a multiple of $\ell(^\varphi a)$. We choose a word
  $b \in  \Gamma^+(c_a), 
 \ell(b)> 2L+1$,  and such that $\bold\Phi(b)$ is synchronizing.
We set
Given these choices we construct the map $\Theta^-_{\langle H \rangle}(X,\widetilde X)$ as follows:
 For $(\delta, a) \in \Omega^-_{\langle  H \rangle}(X)$ let $x\in X$ be such that
 $$
 x_{(-\ell(bc_a d_{(\delta, a)}) , 0]} = bc_a d_{(\delta, a)} , \quad x_{[1, \infty)}= p^{(+, a)},
 $$
 and set
 $$
 \widetilde x = \varphi (x).
 $$
 By the choice of $\widetilde \sigma_a$ and of $ c_a$ there is an index
 $
 i\in [l-\ell(d), L],
 $
 such that
 \begin{align*}
 \widetilde x_{(i , \infty)} = (p^{(-,^\varphi a)}   , p^{(+, ^\varphi a)}  )_{_{(i , \infty)}}
 \quad
  \widetilde x_{[i , \infty)} \neq (p^{(-,^\varphi a)}   , p^{(+, ^\varphi a)}  )_{_{[i , \infty)}}. 
  \tag{6.1}
  \end{align*}
 Set
 $$
 \widetilde a = \widetilde x_{(i, i + \ell(^\varphi a)]}.
 $$
 By the choice of $b$ one has that 
 $$
  \widetilde \delta =[\widetilde x_{(-\ell(d) - L - \ell(b) - 1, i]}]_\approx 
  \in \Delta_{synchro}(  \widetilde X),
 $$
 and by (6.1) 
 $$
 (\widetilde \delta,  \widetilde a ) \in \Pi^-_{\langle  H \rangle}(X).
 $$
The choices for the construction of $\Theta^+_{\langle  H \rangle}(X,\widetilde X)$ and the construction itself are symmetric.

 Conditions  $A(-), A(+)$ and $B(-),B(0),B(+)$ are satisfied by construction.
 \end{proof}
 
For a sofic system $X$ denote by $X^{(n)}$ the subshift of finite type 
that has as its language of admissible words  
set of synchronizing words of $X$ of length $n$. By 
\cite [Proposition 3]{M} one has
$$
\lim_{n \to\infty} \ent (X^{(n)})   =  \ent (X).
$$
 
For the next theorem compare and its proof. 

\begin{theorem}   
Le $X$ and $\widetilde X$ be topologically mixing sofic systems, let 
$\ent (X) >\ent ( \widetilde X)$, and let  $n\in \Bbb  N$ be such that 
$\ent(X^{(n)}) >\ent ( \widetilde X)$.
Let 
\begin{align*}
H \geq n, \tag {6.2}
\end{align*}
and such that (4.2) holds. Let there be given a homomorphism
$$
\varphi_\circ: 
P_{\langle  H \rangle}(X)\to P_{\langle  H \rangle}( \widetilde X),
$$
that is accompanied by a triple 
$(\Theta_{\langle  H \rangle},\Theta_{\langle  H \rangle},\Theta_{\langle  H \rangle})
(X,\widetilde X)$. Then there exists a surjective homomorphism 
$$
\zeta: X\to\widetilde X.
$$
\end{theorem}
\begin{proof}
Denote by $Q(X^{(n)})$ the set of $p \in P(X)$, whose period does not have a divisor among the periods of the Fischer cover of $\widetilde X$. 
From the hypothesis, that
 $\widetilde X$ is topologically mixing, it follows, that its Fischer cover is also topologically mixing. This implies that the set $Q(X^{(n)})$ is finite.
Choose an $m\in \Bbb N$ such that removing the words in
$$
\{p_{[1, mn\pi(p)]}: p \in Q(X^{(n)})  \}
$$
yields an irreducible subshift of finite type $X^{\circ}$ such that $h(X^{\circ})  > h(\widetilde X)$. Without loss of generality one can assume that $m=1$.
By  \cite[Corollary (2.6)]{B} there exist a surjective homomorphism $\chi$ of $X^{\circ}$ onto $\widetilde X$.
 Let $[-l^\circ,l^\circ]$ be a coding window for $\chi$, and let $\chi$ be given by the block map
 $$
\bold{X}: X^\circ_{[-l^\circ,l^\circ]}\to \widetilde  \Sigma. 
$$

We extend the surjective homomorphism $\chi:X^{\circ}\to\widetilde X$ to a homomorphism $$
\zeta: X \to  \widetilde X.
$$
For this we choose  an $ L \in \Bbb N$, such that for 
 $u, u^\prime \in \mathcal L_n(X^{\circ} )$ there exists a word $v \in  \mathcal L(X^{\circ})$
such that $uvu^\prime  \in  \mathcal L(X^{\circ})$,
and also choose an $ \widetilde L\in \Bbb N$, such that for 
$ \widetilde u, \widetilde  u^\prime \in \mathcal L(\widetilde X )$ there exists a word 
$v \in \mathcal L_{\widetilde m}(X) $ such that
$
\widetilde u\widetilde v \widetilde u \in \mathcal L(\widetilde X ).
$

We denote by $\widetilde D^-$($\widetilde D^+$) the smallest 
$\widetilde D\in \Bbb N$, such that we can select for 
$(\widetilde\delta^{(-)},\widetilde a^{(-)} ) \in \Pi^{(-)}(\widetilde X)$
($(\widetilde a^{(+)},\widetilde\delta^{(+)}) \in \Pi^{(+)}(\widetilde X)$) a word
$d(\widetilde\delta^{(-)},\widetilde a^{(-)} )  \in \widetilde\delta^{(-)}$
($d(\widetilde a^{(+)} ,\widetilde\delta^{(+)})  \in \widetilde\delta^{(+)} $)  of length less than or equal to $\widetilde D^-$
($\widetilde D^+(\widetilde X) $), that has a final (initial) symbol, that is different from the final (initial) symbol of $\widetilde a^{(-)}$($\widetilde a^{(+)}$).

We choose a word $u \in \mathcal L(X^{\circ}),\ell(u) > 2l^\circ$, such that 
$\bold{X}(u)$ is synchronizing, and we set $\widetilde u = \bold{X}(u).$

We set
$$
 \widetilde M = \widetilde L + \max \{ \widetilde D^-, \widetilde D^+\}, \quad M = \widetilde M + \ell( \widetilde u ) +
  2l^\circ.
$$

For $v \in \mathcal L_n(X^\circ)$ we choose a $b^-(v)\in  \mathcal L_L(X^\circ)$, such that 
$vb^-(v)u \in  \mathcal L(X^\circ)$. The meaning of $b^-$ is symmetric.
We choose for $\widetilde\delta \in \Delta_{synchro} (X)$  and for 
$
m \geq \widetilde C  +  \widetilde D - \ell(\widetilde d(\widetilde\delta)) 
$
a
$\widetilde b^-(m, \widetilde\delta )  \in \mathcal L_n(\widetilde X)$, such that 
$\widetilde u  \widetilde b^-(m, \widetilde\delta )\widetilde d(\widetilde\delta)   \in 
 \mathcal L(\widetilde X)$.
The meaning of $\widetilde b^+$ is symmetric.

Let there be given $x \in X$. We construct the image of $x$ under $\zeta$, which we denote by 
 $\widetilde z$. 
At hand are a compact-open set $F\subset X$ such that (2.2) and (2.3) hold, and  a  
compact-open set $\widehat F\subset X$,  as in (4.3). Also at hand are time points 
$T(x), x \in\widehat F,$ as in (4.6) and the triples $(a[-],c,a[+])(x), x \in\widehat F$, as in (4.10). Also
the homomorphism $\varphi: X \to\widetilde X$, as constructed in Section 4 can be assumed as given, in particular
the image $\widetilde x$ of $x$ under $\varphi$.

We set
\begin{align*}
&E^-(X) = \{x\in X:x_{[-2M-1, 0)}\in   \mathcal L(X^{\circ}),
x_{(-2H, 0]} \not \in \mathcal L(X^{\circ})\},
\\
&E^+(X) = \{x\in X:x_{( 0.2M +1]}\in   \mathcal L(X^{\circ}),
x_{[0, 2H)} \not \in \mathcal L(X^{\circ})\}.
\end{align*}
As a consequence of (6.2), one has that
$$
x_{(-2H, 0] }\in \mathcal L_{2H}(X)\setminus \mathcal P_{2H}(X), \qquad   x \in E^-(X),
$$
$$
x_{[0,-2H) }\in \mathcal L_{2H}(X)\setminus \mathcal P_{2H}(X), \qquad   x \in E^+(X) ,
$$
and by (2.2) one obtains indices 
$$
j_{F}^{-}(x) = \min \{j \in (-3H, H): S_X^{j}(x) \in F \}, \qquad   x \in E^-(X),
$$
$$
j_{F}^{+}(x) = \max \{j \in (-H, 3H): S_X^{-j}(x) \in F \}, \qquad   x \in E^+(X),
$$
and, as is seen from the structure of $\widehat F$, one has  indices 
$$
j_{\widehat F}^{-}(x) = \min \{ \widehat j \in (j_{F}^{-}(x) - 3H, j^{-}_{F}(x) ): S_X^{\widehat j}(x) \}, 
 \qquad x   \in E^-(X),
$$
$$
 j_{\widehat F}^{+}(x) = \max \{\widehat  j \in ( j_{F}^{+}(x), j_{F}^{+}(x) + 3H ): 
 S_X^{-\widehat j}(x) \},  \qquad x\in E^+(X).
$$
We set
$$
t^-(x) = T(S^{j_{F}^{(-)}}_X (x)),
$$
$$
(a^{(-)}[-],c^{(-)},a^{(-)}[+])(x) = (a[-],c,a[+])(S^{j_{F}^{(-)}}_X (x)),
$$
$$
\widehat t^-(x) = T(S^{ j_{\widehat F}^{-}(x)}_X   (x)),
$$
$$
(\widehat a^{(-)}[-],\widehat c^{(-)},\widehat a^{(-)}[+])(x) = 
(a[-],  c, a[+])(S^{j_{\widehat F}^{-}(x)}_X   (x)), \qquad x \in E^-(X),
$$
$$
t^+(x) = T(S^{j_{F}^{(+)}}_X   (x)),
$$
$$
(a^{(+)}[-],c^{(+)},a^{(+)}[+])(x) = (a[-],c,a[+])(S^{j_{F}^{(+)}}_X (x)),
$$
$$
 \widehat t^+(x) = T(S^{j_{\widehat F}^{+}(x)}_X   (x)),
$$
$$
(\widehat a^{(+)}[-],\widehat c^{(+)},\widehat a^{(+)}[+])(x) = 
(a[-],  c, a[+])(S^{ j_{\widehat F}^{+}(x)}_X   (x)), \qquad x \in E^+(X). 
$$

One obtains pairs $(\delta^{(-)},a^{(-)})(x) \in \Pi^-(X),x\in E^-(X),$ and 
$(a^{(+)},\delta^{(+)})(x) \in \Pi^-(X), x\in E^+(X)$, by 
\begin{align*}
&\delta^{(-)}(x)= [x_{(\widehat t(x) -n,\widehat t(x)]}]_\approx, \quad
a^{(-)}(x) = \widehat a^{(-)}[+](x), \qquad  \  \  \  \  \  \  \  \  x\in E^-(X),
\\
&a^{(+)}(x) = a[+](x), \quad  \delta^{(+)}(x) = [c^{(+)}(x)
x_{(\widehat t^{(+)}(x),
t(x) + n]}]_\approx, \qquad x\in E^+(X).
\end{align*}

For the coding instructions we introduce notation.
We denote the image of the pair $(\delta^{(-)},a^{(-)})(x)$ under 
$\Theta^-_{\langle  H \rangle}(X,\widetilde X)$  by 
$
((\widetilde\delta^{(-)},\widetilde a^{(-)}),\widetilde  t^-)(x),x\in E^-(X),
$
and we denote the image of the pair $(a^{(+)},\delta^{(+)})(x)$ under 
$\Theta^+_{\langle  H \rangle}(X,\widetilde X)$  by 
$
((\widetilde a^{(+)},\widetilde\delta^{(+)}) ,\widetilde  t^+)(x)$, $x\in E^+(X).
$
We denote the image of the triple $(a^{(-)}[-],c^{(-)},a^{(-)}[+])(x)$ under 
$\Theta_{\langle  H \rangle}(X,\widetilde X)$  by 
$$
(\widetilde a^{(-)}[-],\widetilde c^{(-)},\widetilde a^{(-)}[+]), \widetilde t^{(-)} )(x),
\qquad x\in E^-(X),
$$
and we denote the image of the triple $(a^{(+}[-],c^{(+)},a^{(+)}[+])(x)$  under 
$\Theta_{\langle  H \rangle}(X,\widetilde X)$  by 
$$
(\widetilde a^{(+)}[-],\widetilde c^{(+)},\widetilde a^{(+)}[+])(x), \widetilde t^{(+)} )(x),
\qquad x\in E^+(X).
$$
We set
\begin{align*}
&\widehat s^{(-)}(x) =\widehat t^{(-)}(x) + \widetilde t^{-}(x), \
s^{(-)}(x) = t^{(-)}(x) + \widetilde t^{(-)}(x), \qquad x \in E^-(X),
\\
 &s^{(+)}(x) = t^{(-)}(x) + \widetilde t^{(+)}(x), \
\widehat s^{(+)}(x) =\widehat t^{(+)}(x) + \widetilde t^{+}(x), \qquad x \in E^+(X),
\end{align*}

\medskip
\begin{align*}
 q^{(-)} (x) = \lfloor \widehat s^{(+)}(x)-s^{(+)}(x)   /
\gamma(\widetilde a^{(+)}(x)) \ell(\widetilde a^{(+)}(x))\rfloor, 
&
\\r^{(-)} (x) = \ s^{(-)}(x) - \widehat s^{(-)}(x) - 
\ell(\widetilde c^{(-)}(x))
-q^{(-)} (x) \
&\gamma(\widetilde a^{(-)}(x))\ell(\widetilde a^{(-)}(x)), 
\\
&
\qquad \qquad \   \   
 x \in E^-(X),
\\
\\
q^{(+)} (x) = \lfloor \widehat s^{(+)}(x)-s^{(+)}(x)   /
\gamma(\widetilde a^{(+)}(x)) \ell(\widetilde a^{(+)}(x))\rfloor, \ &
\\
 r^{(+)} (x) = \ \widehat s^{(+)}(x) -s^{(+)}(x)  \  -q^{(+)} (x) 
\gamma(\widetilde a^{(+)}(x))\ell(\widetilde a&^{(+)}(x)),
 \quad x \in E^+(X).
\end{align*}

\medskip

\noindent
We construct $\widetilde z$ partially. We set 
\begin{align*}
&\widetilde z_{[- M - l_\circ,-\widetilde M]}=\bold X(x_{[- M - l_\circ, - M + l_\circ]}
 b^-(x_{[- M - l_\circ, - M + l_\circ]},u)u), 
\\
&\widetilde z_{( -\widetilde M,\widehat s^{(-)}(x)]} =
\widetilde b^-(\widehat s^{(-)}(x) +\widetilde M + 6H -
\ell(\widetilde d((\widetilde\delta^{(-)}, \widetilde a^{(-)})(x))
\widetilde d((\widetilde\delta^{(-)}, \widetilde a^{(-)})(x)),
\\
&\widetilde z_{( \widehat s^{(-)}(x), s^{(-)}(x)]} = 
\widetilde a^{(-)}(x)^{q^{(-)}(x)\gamma(\widetilde a^{(-)}(x))}
\widetilde a^{(-)}(x)^{\gamma(\widetilde a^{(-)}(x))}_{[1, r^{(-)} (x)]},\quad x \in E^-,
\end{align*}
and for $x \in X$, such that
$$
x_{[-M,M]}\in \mathcal L(x^\circ),
$$
we set
$$
\widetilde z_0 = \bold {X}(x_{[-l_\circ,l_\circ]}),
$$
and we set
\begin{align*}
&\widetilde z_{(  s^{(+)}(x), \widehat s^{(+)}(x)]} = 
\widetilde a^{(+)}(x)^{q^{(+)}(x)\gamma(\widetilde a^{(+)}(x))}
\widetilde a^{(+)}(x)^{\gamma(\widetilde a^{(+)}(x))}_{[1, r^{(+)} (x)]},
\\
&\widetilde z_{[ \widehat s^{(+)}(x), \widetilde M)} =
\widetilde d((\widetilde a^{(+)},\widetilde\delta^{(+)})(x))\widetilde b^+(\widehat s^{(+)}(x) + \widetilde M  + 6H -  \ell(\widetilde d((\widetilde a^{(+)},\widetilde\delta^{(+)})(x)))), 
\\
&\widetilde z_{[\widetilde M, M + l_\circ]} =
\bold X(ub^+(u,x_{[M - l_\circ, M + l_\circ]})
x_{[M - l_\circ, M + l_\circ]}), \quad x \in E^+.
\end{align*}

\medskip
\noindent
One completes the construction of $\widetilde z$ by copying the remaining blocks from 
$\widetilde x$.
For $ x \in E^-(X)$ such that
$$
i^+_{E^+(X)}(S_X^{-n}(x)) < \infty,
$$
and such that 
$$
j^-_F(x) < i^+_{E^+(X)}(S_X^{-n}(x)) + j^+_F(S_X^{- i^+_{E^+(X)}(S_X^{-n}x) -n}(x)),
$$
one sets
$$
\widetilde z_{(s^-(x), s^+(S_X^{- i^+_{E^+(X)}(S_X^{-n}x) -n}(x)))} = 
\widetilde x_{(s^-(x), s^+(S_X^{- i^+_{E^+(X)}(S_X^{-n}x) -n}(x)))}.
$$
For $ x \in E^-(X)$ such that
$$
i^+_{E^+(X)}(S_X^{-n}(x)) = \infty,
$$
 one sets
 $$
 \widetilde z_{(s^-(x), \infty)} = 
\widetilde x_{(s^-(x), \infty)},
 $$
 and for 
 $ x \in E^+(X)$ such that
$$
i^-_{E^-(X)}(S_X^{n}(x)) = \infty,
$$
 one sets
 $$
 \widetilde z_{(\infty, s^+(x))} = 
\widetilde x_{(\infty, s^+(x))}.
$$
Conditions $(A)$, $(A(-))$, $(A(+))$, and $(B)$, $(B(-))$,  $(B(0))$, $(B(+))$ assure that $\widetilde z$ is well defined and is in $\widetilde X$. One also checks that the map 
$x \to \widetilde z$ has a coding window.
\end{proof}
The considerations on decidability of Section 4 extend to the case of surjective homomorphisms.

\bigskip

\par\noindent Wolfgang Krieger
\par\noindent Institute for Applied Mathematics, 
\par\noindent  University of Heidelberg,
\par\noindent Im Neuenheimer Feld 205, 
 \par\noindent 69120 Heidelberg,
 \par\noindent Germany
\par\noindent krieger@math.uni-heidelberg.de


\begin{thebibliography}{9999)}
 
\bibitem[AGW]{AGW}{\sc R. ~Adler, L. W. ~Goodwyn and B. ~Weiss},
{\it  Equivalence of topological Markov shifts}, 
Israel J. of Math.
{\bf 27}
(1977),
 49 -- 63
 
 \bibitem[AM]{AM}
{\sc  R.~Adler  and B.~Marcus},
{\it  Topological entropy and equivalence of dynamical systems},
Memoirs of the American Mathematical Society
{\bf  219}
(1979)

 \bibitem[AKM]{AKM}{\sc R. ~Adler, B. ~Kitchens and B.~ Marcus},
{\it  Almost topological classification of finite-to-one factor maps between shifts of finite type}, 
Ergod Th. \& Dynam. Sys.
{\bf 5}
(1985),
 485 -- 500
  
\bibitem[B]{B}{\sc M. ~Boyle},
{\it  Lower entropy factors of sofic systems}, 
Ergod Th. \& Dynam. Sys.
{\bf 4}
(1984),
 541 -- 557

\bibitem[CP]{CP}
{\sc  E.~Coven and M.~Paul},
{\it Finite procedures for sofic systems },
Monatsh. Math.
{\bf  83}
(1977),
205 -- 278

\bibitem [F]{F}
{\sc R.~Fischer},
{\it  Sofic systems and graphs},
Monatsh. Math. 
{\bf 80}
(1975),
179 -- 186

\bibitem[HU]{HU}{\sc J.E.~Hopcroft and J.D.~Ullman},
{\it Introduction to Automata Theory, Languages, and Computation},
 Addison-Wesley, Reading
(2001)

\bibitem[Ki]{Ki}{\sc B.~P.~Kitchens},
{\it Symbolic dynamics}, Springer, Berlin, Heidelberg, New York
(1998)
 
 \bibitem [Kr]{Kr}
{\sc W.~Krieger},
{\it On the subsystems of topological Markov chains},
Ergod Th. \& Dynam. Sys.
{\bf  2}
(1982),
195 -- 202

\bibitem[LM]{LM}{\sc D.~Lind and B.~Marcus},
{\it An introduction to symbolic dynamics and coding},
 Cambridge University Press, Cambridge
(1995)

\bibitem[M]{M}{\sc B.~Marcus},
{\it Sofic Systems and Encoding Data},
IEEE Tansactions on Information Theory
{\bf  31}
(1985),
366  -- 377

\bibitem[N]{N}{\sc J. A.~Nielsen},
{\it Morphisms between sofic shift spaces},
Thesis, Department of Mathematics, Aarhus University,
(2010)

 \bibitem [P]{P}
{\sc W.~Parry},
{\it Intrinsic Markov chains},
Trans. Amer. Math. Soc.
{\bf  112}
(1964),
55 -- 66

\bibitem [PT]{PT}{\sc W.~Parry and S.~Tuncel},
{\it Classification Problems in Ergodic Theory}, Cambridge University Press
(1982)

\bibitem [T]{T}
{\sc K.~Thomsen},
{\it On the structure of the sofic shift space},
Trans. Amer. Math. Soc.
{\bf  365}
(2004),
3557 -- 3619

\bibitem [W]{W}
{\sc B.~Weiss},
{\it Subshifts of finite type and sofic systems},
Monatsh. Math.
{\bf  77}
(1973),
462 -- 474

\end{thebibliography}
\end{document}